\documentclass{amsproc}%
\usepackage{graphicx}
\usepackage{amscd}
\usepackage{amsmath}
\usepackage{amsfonts}
\usepackage{amssymb}%
\theoremstyle{plain}

\numberwithin{equation}{section}

\begin{document}
\title[On the left and right matrix coefficients]{On the left and right coefficients of the general Cayley-Hamilton identities
for an $n\times n$ matrix}
\author{Szilvia Homolya}
\address{Institute of Mathematics, University of Miskolc \\
3515 Miskolc-Egyetemv\'{a}ros, Hungary}
\email{szilvia.homolya@uni-miskolc.hu}
\author{Jen\H{o} Szigeti}
\address{Institute of Mathematics, University of Miskolc \\
3515 Miskolc-Egyetemv\'{a}ros, Hungary}
\email{jeno.szigeti@uni-miskolc.hu}
\thanks{The second named author was partially supported by the National Research,
Development and Innovation Office of Hungary (NKFIH) K138828}

\begin{abstract}
An $n\times n$ matrix $A\in\mathrm{M}_{n}(R)$ over an arbitrary unitary ring
$R$ satisfies the following invariant left and right Cayley-Hamilton
identities
\[
(\mu_{0}I_{n}+C(0))+(\mu_{1}I_{n}+C(1))A+\cdots+(\mu_{n-1}I_{n}+C(n-1))A^{n-1}%
+(n!I_{n}+C(n))A^{n}=0,
\]%
\[
(\mu_{0}I_{n}+D(0))+A(\mu_{1}I_{n}+D(1))+\cdots+A^{n-1}(\mu_{n-1}%
I_{n}+D(n-1))+A^{n}(n!I_{n}+D(n))=0,
\]
where $C(i),D(i)\in\mathrm{M}_{n}([R,R]),0\leq i\leq n$ are matrix
coefficients with commutator sum entries and $\mu_{0}+\mu_{1}t+\cdots
+\mu_{n-1}t^{n-1}+n!t^{n}\in R[t]$\ is the so-called symmetric characteristic
polynomial of $A$. If $R=R_{0}\oplus R_{1}$ is a $\mathbb{Z}_{2}$-grading
similar to the case of Grassmann algebras, then we prove that
$C(i)-D(i)-AC(i+1)+D(i+1)A=-2\mu_{i+1}^{(1)}A_{1}$ for all $-1\leq i\leq n$,
where $A_{1}\in\mathrm{M}_{n}(R_{1})$ and $\mu_{i+1}^{(1)}\in R_{1}$ are the
odd components of $A$ and $\mu_{i+1}$ respectively.\textit{ }

\end{abstract}
\maketitle

\noindent1. INTRODUCTION\ AND\ BASIC\ DEFINITIONS

\bigskip

\noindent The natural symmetrization of the determinant and the adjoint matrix
can serve as a starting point of a symmetric variant of the classical
determinant theory for square matrices over an arbitrary ring. The most
important feature of this approach is that it works for square matrices over
an arbitrary ring and especially over a Lie nilpotent ring.

\noindent The Grassmann (exterior) algebra%
\[
E=K\left\langle v_{1},v_{2},...,v_{i},...\mid v_{i}v_{j}+v_{j}v_{i}=0\text{
for all }1\leq i\leq j\right\rangle
\]
(over a field $K$) generated by the infinite sequence of anticommutative
indeterminates $(v_{i})_{i\geq1}$ is our basic example, which is Lie nilpotent
of index $2$.

\noindent The natural $\mathbb{Z}_{2}$-grading is $E=E_{0}\oplus E_{1}$ with
$E_{0}$ being the $K$-subspace (subalgebra) generated by the monomials of even
length and $E_{1}$ being the $K$-subspace generated by the monomials of odd
length. Any element $g\in E$ can be uniquely written as $g=g_{0}+g_{1}$, where
$g_{0}\in E_{0}$ is the even and $g_{1}\in E_{1}$ is the odd component of $g$.
Notice that $g_{0}\in E_{0}=\mathrm{Z}(E)$ is central and $(g_{1})^{2}=0$.

\noindent In case of $\mathrm{char}(K)=0$, Kemer's pioneering work on the
T-ideals of associative algebras revealed the importance of the identities
satisfied by the full $n\times n$ matrix algebra $\mathrm{M}_{n}(E)$ and by
its subalgebra of $(n,t)$ supermatrices $\mathrm{M}_{n,t}(E)$.

\noindent Let $K\left\langle x_{1},x_{2},\ldots,x_{i},\ldots\right\rangle $
denote the polynomial $K$-algebra generated by the infinite sequence
$x_{1},x_{2},\ldots,x_{i},\ldots$\ of non-commuting indeterminates. The prime
T-ideals of this (free associative $K$-)algebra\ are exactly the T-ideals of
the identities satisfied by $\mathrm{M}_{n}(K)$ for $n\geq1$ (already proved
by Amitsur in [A]). The T-prime T-ideals are the prime T-ideals plus the
T-ideals of the identities of $\mathrm{M}_{n}(E)$ for $n\geq1$ and of
$\mathrm{M}_{n,t}(E)$ for $n-1\geq t\geq1$ (see [K]). Another remarkable
result is that for a sufficiently large $n\geq1$, any T-ideal contains the
T-ideal of the identities satisfied by $\mathrm{M}_{n}(E)$ (see [K]). The
above mentioned algebras $\mathrm{M}_{n}(E)$ and $\mathrm{M}_{n,t}(E)$ served
as the main motivation for the development of the symmetric determinant.

\noindent For an $n\times n$ matrix matrix $A$ over a Lie nilpotent ring $R$
of index $k$, the mentioned symmetric determinant can be used to exhibit an
invariant Cayley-Hamilton identity of the form%
\[
\lambda_{0}^{(k)}I_{n}+\lambda_{1}^{(k)}A+\cdots+\lambda_{n^{k}-1}%
^{(k)}A^{n^{k}-1}+\lambda_{n^{k}}^{(k)}A^{n^{k}}=0
\]
with left scalar coefficients $\lambda_{i}^{(k)}\in R$, $0\leq i\leq n^{k}$
and $\lambda_{n^{k}}^{(k)}=n\left\{  (n-1)!\right\}  ^{1+n+\cdots+n^{k-1}}$.
If $k=2$, then the following two important consequences can be derived:
$\mathrm{M}_{n}(E)$ is integral of degree $2n^{2}$ and $\mathrm{M}_{n,d}(E)$
is integral of degree $n^{2}$ over the even component $E_{0}$ of $E$ (see [S1,
S2]). The above Cayley-Hamilton identity also can be used to prove the left
(and right) integrality of a Lie nilpotent ring over the fixed ring of any of
its finite automorphism groups (see [S4]). As we can see the details in
section 3, for an $n\times n$ matrix over an arbitrary ring, the use of the
symmetric characteristic polynomial yields a general left (and right)
invariant Cayley-Hamilton identity (of degree $n$)\ with certain special
matrix coefficients (see [S3]).

\noindent In what follows an algebra $R$\ means a not necessarily commutative
unitary algebra over a field (commutative unitary ring), the notation for the
full $n\times n$ matrix algebra over $R$ is $\mathrm{M}_{n}(R)$ and $I_{n}%
\in\mathrm{M}_{n}(R)$ denotes the identity. The matrix $E_{i,j}\in
\mathrm{M}_{n}(R)$ has $1$ in the $(i,j)$ position and zeros in all other
positions. Let $\mathrm{S}_{n}=\mathrm{Sym}(\{1,\ldots,n\})$ denote the
symmetric group of all permutations of the set $\{1,2,\ldots,n\}$. If%
\[
A=[a_{i,j}]=\underset{1\leq i,j\leq n}{\sum}a_{i,j}E_{i,j}%
\]
is an $n\times n$ matrix over $R$, then%
\[
\mathrm{sdet}(A)=\underset{\tau,\rho\in\mathrm{S}_{n}}{\sum}\mathrm{sgn}%
(\rho)a_{\tau(1),\rho(\tau(1))}\cdots a_{\tau(t),\rho(\tau(t))}\cdots
a_{\tau(n),\rho(\tau(n))}%
\]%
\[
=\underset{\alpha,\beta\in\mathrm{S}_{n}}{\sum}\mathrm{sgn}(\alpha
)\mathrm{sgn}(\beta)a_{\alpha(1),\beta(1)}\cdots a_{\alpha(t),\beta(t)}\cdots
a_{\alpha(n),\beta(n)}%
\]
is called the symmetric determinant of $A$. If $R$\ is commutative, then
$\mathrm{sdet}(A)=n!\det(A)$, where $\det(A)$ denotes the ordinary (classical)
determinant of $A$.

\noindent The $(r,s)$-entry of the symmetric adjoint matrix $A^{\ast}%
=[a_{r,s}^{\ast}]$ of $A$ is defined as follows:%
\[
a_{r,s}^{\ast}=\underset{\tau,\rho}{\sum}\mathrm{sgn}(\rho)a_{\tau
(1),\rho(\tau(1))}\cdots a_{\tau(s-1),\rho(\tau(s-1))}a_{\tau(s+1),\rho
(\tau(s+1))}\cdots a_{\tau(n),\rho(\tau(n))}%
\]%
\[
=\underset{\alpha,\beta}{\sum}\mathrm{sgn}(\alpha)\mathrm{sgn}(\beta
)a_{\alpha(1),\beta(1)}\cdots a_{\alpha(s-1),\beta(s-1)}a_{\alpha
(s+1),\beta(s+1)}\cdots a_{\alpha(n),\beta(n)}\text{ },
\]
where the first sum is taken over all $\tau,\rho\in\mathrm{S}_{n}$ with
$\tau(s)=s$ and $\rho(s)=r$, while the second sum is taken over all
$\alpha,\beta\in\mathrm{S}_{n}$ with $\alpha(s)=s$ and $\beta(s)=r$. Clearly,
the $(r,s)$ entry of $A^{\ast}$ is exactly the signed symmetric determinant
$(-1)^{r+s}\mathrm{sdet}(A_{s,r})$\ of the $(n-1)\times(n-1)$\ minor $A_{s,r}%
$\ of $A$ arising from the deletion of the $s$-th row and ther $r$-th column
of $A$. An immediate corollary is that for a commutative base ring $R$\ we
have $A^{\ast}=(n-1)!\mathrm{adj}(A)$, where $\mathrm{adj}(A)$ denotes the
ordinary adjoint of $A$.

\noindent If $R$ is commutative, then $\mathrm{tr}(AB)=\mathrm{tr}(BA)$\ for
all $A,B\in\mathrm{M}_{n}(R)$. In spite of the fact that this well known trace
identity is no longer valid for matrices over a non-commutative ring, we still
have the following.

\bigskip

\noindent\textbf{Theorem} (see [SvW])\textbf{.}\textit{ The traces of the
product matrices }$A^{\ast}A$\textit{ and }$AA^{\ast}$\textit{ are both equal
to the symmetric determinant of }$A$\textit{:}%
\[
\mathrm{tr}(AA^{\ast})=\mathrm{sdet}(A)=\mathrm{tr}(A^{\ast}A).
\]

\bigskip

\noindent The next theorem of Domokos is about the invariant property of the
symmetric adjoint and determinant.

\bigskip

\noindent\textbf{Theorem} (see [D])\textbf{. }\textit{Let }$R$\textit{ be an
algebra over a field }$K$\textit{ of characteristic zero. If }$A\in
\mathrm{M}_{n}(R)$\textit{ and }$P\in\mathrm{GL}_{n}(K)$\textit{ is
invertible, then }$(PAP^{-1})^{\ast}=PA^{\ast}P^{-1}$\textit{ and
}$\mathrm{sdet}(PAP^{-1})=\mathrm{sdet}(A)$\textit{.}

\bigskip

\noindent The symmetric (as well as the Lie nilpotent) version of classical
determinant theory rests on the following result.

\bigskip

\noindent\textbf{Theorem} (see [S1])\textbf{.}\textit{ For a matrix }%
$A\in\mathrm{M}_{n}(R)$\textit{ the entries of the left and the right
commutator parts}%
\[
C=nA^{\ast}A-\mathrm{tr}(A^{\ast}A)I_{n}=nA^{\ast}A-\mathrm{sdet}%
(A)I_{n}\mathit{\ }%
\]%
\[
D=nAA^{\ast}-\mathrm{tr}(AA^{\ast})I_{n}=nAA^{\ast}-\mathrm{sdet}(A)I_{n}%
\]
\textit{of }$A$\textit{\ are in the additive commutator subgroup }%
$[R,R]$\textit{ of }$R$\textit{ generated by the elements of the form
}$[u,v]=uv-vu$\textit{ with }$u,v\in R$\textit{. We also have }$\mathrm{tr}%
(C)=\mathrm{tr}(D)=0$\textit{.}

\bigskip

\noindent\textbf{Remarks.} We note that each entry of $C$ and $D$\ is a sum of
commutators of the form $\pm\lbrack a^{\prime},a^{\prime\prime}]$, where
$a^{\prime}$ and $a^{\prime\prime}$ are products of certain entries of $A$.
Using%
\[
A=\left[
\begin{array}
[c]{cc}%
a & b\\
c & d
\end{array}
\right]  \text{ and }A^{\ast}=\left[
\begin{array}
[c]{cc}%
d & -b\\
-c & a
\end{array}
\right]  ,
\]
the $2\times2$ case is illustrated in the following calculations.%
\[
C=2A^{\ast}A-\mathrm{tr}(A^{\ast}A)I_{2}=2A^{\ast}A-\mathrm{sdet}%
(A)I_{2}=\left[
\begin{array}
[c]{cc}%
-[a,d]+[c,b] & 2[d,b]\\
2[a,c] & [a,d]-[c,b]
\end{array}
\right]  ,
\]%
\[
D=2AA^{\ast}-\mathrm{tr}(AA^{\ast})I_{2}=2AA^{\ast}-\mathrm{sdet}%
(A)I_{2}=\left[
\begin{array}
[c]{cc}%
\lbrack a,d]+[c,b] & 2[b,a]\\
2[c,d] & -[a,d]-[c,b]
\end{array}
\right]  .
\]

\newpage

\noindent2. SOME IDENTITIES\ INVOLVING\ COMMUTATORS AND\ A\ CONNECTION

\noindent\ \ \ BETWEEN\ THE COMMUTATOR\ PARTS $C$ AND $D$

\bigskip

\noindent For a sequence $x_{1},x_{2},\ldots,x_{m}$ of elements in $R$ we use
the notation $[x_{1},x_{2},\ldots,x_{m}]$\ for the left normed commutator
(Lie-)product:%
\[
\lbrack x_{1}]=x_{1}\text{ and }[x_{1},x_{2},\ldots,x_{m}]=[\ldots
\lbrack\lbrack x_{1},x_{2}],x_{3}],\ldots,x_{m}].
\]
A ring $R$ is called Lie nilpotent of index $k\geq1$ if $[x_{1},x_{2}%
,\ldots,x_{k+1}]=0$ is an identity on $R$. A ring $R$ satisfies the Engel
condition (Engelian) of index $k$, if%
\[
\lbrack x,\underset{k}{\underbrace{y,\ldots,y}}]=0
\]
is an identity on $R$. The next (weak) version of Latyshev's theorem (see [L])
reveals an interesting consequence of the Engelian property.

\bigskip

\noindent\textbf{2.1. Theorem.}\textit{ If }$R$\textit{ is an Engelian algebra
(over an infinite field) of index }$k$\textit{, generated by }$m$\textit{
elements, then there exists an integer }$d=d(k,m)$\textit{ such that }%
$R$\textit{\ satisfies the polynomial identity }$[x_{1},y_{1}][x_{2}%
,y_{2}]\cdots\lbrack x_{d},y_{d}]=0$\textit{ (notice that}

\noindent$\lbrack x_{1},y_{1}]z_{1}[x_{2},y_{2}]z_{2}\cdots z_{d-1}%
[x_{d},y_{d}]=0$\textit{ is an easy consequence of it). In the original
version }$R$\textit{\ satisfies a so-called non-matrix polynomial identity.}

\bigskip

\noindent\textbf{Remarks.} The identity $[x_{1},y_{1}]z_{1}[x_{2},y_{2}%
]z_{2}\cdots z_{d-1}[x_{d},y_{d}]=0$ on $R$\ is obviously equivalent to the
fact that the commutator ideal $R[R,R]R\vartriangleleft R$ is nilpotent:
$(R[R,R]R)^{d}=\{0\}$. We note that the Engel condition of index $k$ plus
finite generation imply Lie nilpotency of index $l\geq k$ (see [RW]).

\bigskip

\noindent The already mentioned (infintely generated)\ exterior algebra $E$
over a field $K$ (with $2\neq0$) is Lie nilpotent of index $2$ and%
\[
\lbrack v_{1},v_{2}]\cdot\lbrack v_{3},v_{4}]\cdot\ \cdots\ \cdot\lbrack
v_{2d-1},v_{2d}]=2^{d}v_{1}v_{2}\cdots v_{2d}\neq0
\]
shows that Lie nilpotency of index $2$ alone does not imply $[x_{1}%
,y_{1}][x_{2},y_{2}]\cdots\lbrack x_{d},y_{d}]=0$ for any $d$. On the other
hand, the following implication holds.

\bigskip

\noindent\textbf{2.2. Theorem.}\textit{ If }$R$\textit{\ is Engelian of index
}$k+1$\textit{ and }$\frac{1}{d}\in R$\textit{ for }$d=\binom{2k}{k}$\textit{,
then }%
\[
\lbrack y_{1},\underset{k}{\underbrace{x,\ldots,x}}][y_{2},\underset
{k}{\underbrace{x,\ldots,x}}]=0
\]
\textit{is an identity on }$R$\textit{.}

\bigskip

\noindent\textbf{Proof.} The following identity can be found in [HSvWZ]: for
all $r_{1},r_{2},x_{1},\ldots,x_{m}~\in~R$%
\[
\lbrack r_{1}r_{2},x_{1},\ldots,x_{m}]=\underset{j_{1}<j_{2}<\cdots<j_{m-q}%
}{\underset{1\leq i_{1}<i_{2}<\cdots<i_{q}\leq m}{%
{\displaystyle\sum}
}}[r_{1},x_{i_{1}},\ldots,x_{i_{q}}]\cdot\lbrack r_{2},x_{j_{1}}%
,\ldots,x_{j_{m-q}}]\text{ \ \ }(\ast)
\]
holds, where the sum is taken over all strictly increasing sequences

\noindent$1\leq i_{1}<i_{2}<\cdots<i_{q}\leq m$ with $0\leq q\leq m$ and%
\[
\{j_{1},j_{2},\ldots,j_{m-q}\}=\{1,2,\ldots,m\}\smallsetminus\{i_{1}%
,i_{2},\ldots,i_{q}\}\text{ and }1\leq j_{1}<j_{2}<\cdots<j_{m-q}\leq m.
\]
The empty and the full sequences ($q=0$ or $q=m$) are allowed with
$[r_{1},\varnothing]=r_{1}$ and $[r_{2},\varnothing]=r_{2}$. Take $m=2k$, then
the substitution%
\[
x_{1}=\cdots=x_{2k}=x
\]
in $(\ast)$\ immediately gives that%
\[
\lbrack r_{1}r_{2},\underset{2k}{\underbrace{x,\ldots,x}}]=\underset
{j_{1}<j_{2}<\cdots<j_{2k-q}}{\underset{1\leq i_{1}<i_{2}<\cdots<i_{q}\leq2k}{%
{\displaystyle\sum}
}}[r_{1},\underset{q}{\underbrace{x,\ldots,x}}]\cdot\lbrack r_{2}%
,\underset{2k-q}{\underbrace{x,\ldots,x}}]\text{ \ \ }(\ast\ast),
\]
where the sum is taken over all strictly increasing sequences

\noindent$1\leq i_{1}<i_{2}<\cdots<i_{q}\leq2k$ with $0\leq q\leq2k$ and%
\[
\{j_{1},j_{2},\ldots,j_{2k-q}\}=\{1,2,\ldots,2k\}\smallsetminus\{i_{1}%
,i_{2},\ldots,i_{q}\}\text{, }1\leq j_{1}<j_{2}<\cdots<j_{2k-q}\leq2k.
\]
The empty and the full sequences ($q=0$ or $q=2k$) are allowed with
$[r_{1},\varnothing]=r_{1}$ and $[r_{2},\varnothing]=r_{2}$. If $R$\ is
Engelian of index $k+1$, then $q\geq k+1$ implies $[r_{1},\underset
{q}{\underbrace{x,\ldots,x}}]=0$ and $2k-q\geq k+1$ implies that
$[r_{2},\underset{2k-q}{\underbrace{x,\ldots,x}}]=0$. Thus the only non-zero
summands in $(\ast\ast)$ are the products $[r_{1},\underset{k}{\underbrace
{x,\ldots,x}}]\cdot\lbrack r_{2},\underset{k}{\underbrace{x,\ldots,x}}]$. The
number of such summands is $d=\binom{2k}{k}$, whence%
\[
0=[r_{1}r_{2},\underset{2k}{\underbrace{x,\ldots,x}}]=d[r_{1},\underset
{k}{\underbrace{x,\ldots,x}}]\cdot\lbrack r_{2},\underset{k}{\underbrace
{x,\ldots,x}}]
\]
follows. In view of $\frac{1}{d}\in R$ our proof is done. $\square$

\bigskip

\noindent\textbf{Remark.} The algebra $\mathrm{U}_{2}(K)$\ of upper triangular
$2\times2$\ matrices is non Engelian of any index $k$, but obviously satisfies
$[x_{1},y_{1}][x_{2},y_{2}]=0$. Clearly,

\noindent$\lbrack\lbrack\lbrack\ldots\lbrack\lbrack E_{1,2},E_{1,1}%
],E_{1,1}],\ldots],E_{1,1}],E_{1,1}]=\pm E_{1,2}\neq0$.

\bigskip

\noindent We mention the following natural generalization of $(\ast)$ in the
above proof of Theorem 2.2.%

\[
\lbrack r_{1}r_{2}\cdots r_{s},x_{1},\ldots,x_{m}]=\underset{(H_{1}%
,H_{2},\ldots,H_{s})}{%
{\displaystyle\sum}
}[r_{1},H_{1}]\cdot\lbrack r_{2},H_{2}]\cdot\cdots\cdot\lbrack r_{s},H_{s}]
\]
holds for all $r_{1},r_{2},\ldots,r_{s},x_{1},\ldots,x_{m}\in R$, where the
sum is taken over all ordered $s$-tuples $(H_{1},H_{2},\ldots,H_{s})$ such
that%
\[
H_{1}\cup H_{2}\cup\cdots\cup H_{s}=\{1,2,\ldots,m\}
\]
is a pair-wise disjoint union ($H_{i}\cap H_{j}=\varnothing$ for all $1\leq
i<j\leq s$). If $r\in R$ and $H\subseteq\{1,2,\ldots,m\}$, then
\[
\lbrack r,H]=[r,x_{h_{1}},\ldots,x_{h_{q}}]\text{ and }[r,\varnothing]=r,
\]
where $H=\{h_{1},h_{2},\ldots,h_{q}\}$ with $1\leq h_{1}<h_{2}<\cdots
<h_{q}\leq m$. It can be proved by using $(\ast)$\ and applying a
straightforward induction with respect to $s$.

\bigskip

\noindent Our first result about the commutator parts is the following.

\bigskip

\noindent\textbf{2.3. Theorem.}\textit{ If }$A\in\mathrm{M}_{n}(R)$\textit{ is
an }$n\times n$\textit{ matrix, then}%
\[
AC-DA=\lambda A-A\lambda=[\lambda I_{n},A]=[[\lambda,a_{i,j}]],
\]
\textit{where }$\lambda=\mathrm{sdet}(A)$\textit{ is the symmetric determinant
and }$C,D\in\mathrm{M}_{n}([R,R])$\textit{ are the left and right commutator
parts of }$A$\textit{. If the base ring }$\frac{1}{2}\in R$\textit{ satisfyies
the Engel identity }$[x,y,y]=0$\textit{ of index }$2$\textit{, then the
product of any two (not necessarily different) entries of the matrix }%
$AC-DA$\textit{ is zero.}

\bigskip

\noindent\textbf{Proof.} Using the notation $\lambda=\mathrm{tr}(AA^{\ast
})=\mathrm{tr}(A^{\ast}A)=\mathrm{sdet}(A)$ and the definitions%
\[
C=nA^{\ast}A-\mathrm{tr}(A^{\ast}A)I_{n}=nA^{\ast}A-\mathrm{sdet}%
(A)I_{n}\mathit{\ }%
\]%
\[
D=nAA^{\ast}-\mathrm{tr}(AA^{\ast})I_{n}=nAA^{\ast}-\mathrm{sdet}(A)I_{n}%
\]
of the commutator parts, we obtain that%
\[
nAA^{\ast}A=A(nA^{\ast}A)=A(C+\lambda I_{n})=AC+A\lambda.
\]
and%
\[
nAA^{\ast}A=(nAA^{\ast})A=(D+\lambda I_{n})A=DA+\lambda A
\]
In consequence, we have%
\[
AC+A\lambda=DA+\lambda A,
\]
whence%
\[
AC-DA=\lambda A-A\lambda=[\lambda a_{i,j}-a_{i,j}\lambda]=[[\lambda,a_{i,j}]]
\]
follows. Thus each entry of $AC-DA$ is a commutator of the form $[\lambda
,a_{i,j}]$. If $R$ satisfies $[x,y,y]=0$, then $[x,r_{1}][x,r_{2}]=0$ is an
identity on $R$ (by Theorem 2.2) and $[\lambda,a_{i_{1},j_{1}}][\lambda
,a_{i_{2},j_{2}}]=0$ follows for all $1\leq i_{1},j_{1},i_{2},j_{2}\leq n$.
$\square$

\bigskip

\noindent An immediate consequence is the following.

\bigskip

\noindent\textbf{Corollary.}\textit{ If }$A\in\mathrm{M}_{n}(R)$\textit{ is an
}$n\times n$\textit{ matrix over a base ring }$\frac{1}{2}\in R$\textit{
satisfying the Engel identity }$[x,y,y]=0$\textit{ of index }$2$\textit{, then
we have }$(AC-DA)^{2}=0$\textit{ for the left and right commutator parts
}$C,D\in\mathrm{M}_{n}([R,R])$\textit{ of }$A$.

\bigskip

\noindent In case of a higher index Engelian property we can present a weaker result.

\bigskip

\noindent\textbf{2.4. Theorem.}\textit{ If }$A\in\mathrm{M}_{n}(R)$\textit{ is
an }$n\times n$\textit{ matrix over an Engelian algebra }$R$\textit{\ (over an
infinite field) of index }$k$\textit{, then for some integer }$d\geq1$,
\textit{the product of any }$d$\textit{ (not necessarily different) entries of
the matrix }$AC-DA$\textit{ is zero (}$C,D\in\mathrm{M}_{n}(R)$\textit{ are
the left and right commutator parts of }$A$\textit{) and }$(AC-DA)^{d}%
=0$\textit{.}

\bigskip

\noindent\textbf{Proof.} Let $S$\ denote the subalgebra of $R$ generated by
the entries $a_{i,j}$ of $A$. Since $S$\ is finitely generated and also
Engelian of index $k$, the application of Latyshev's Theorem 2.1 gives that
$[x_{1},y_{1}][x_{2},y_{2}]\cdots\lbrack x_{d},y_{d}]=0$ (as well as
$[x,y_{1}][x,y_{2}]\cdots\lbrack x,y_{d}]=0$) is an identity on\ $S$\ for some
$d\geq1$. In view of

\noindent$A,C,D\in\mathrm{M}_{n}(S)$ and $\lambda\in S$, our claim is an
immediate consequence of $AC-DA=[[\lambda,a_{i,j}]]$ in Theorem 2.3. $\square$

\bigskip

\noindent A Grassmann like $\mathbb{Z}_{2}$-grading (GL-grading) of a
$K$-algebra $R$ is a $2$-tuple $(R_{0},R_{1})$, where $R_{0},$ and $R_{1}$ are
$K$-subspaces of $R$ such that $R=R_{0}\oplus R_{1}$ is a direct sum and
$R_{i}R_{j}\subseteq R_{i+j}$ for all $i,j\in\{0,1\}$, where $i+j$ is taken in
$\{0,1\}$\ modulo $2$, $R_{0}\subseteq\mathrm{Z}(R)$ is central and
$(r_{1})^{2}=0$ for all $r_{1}\in R_{1}$. If $r=r_{0}+r_{1}$ and
$s=s_{0}+s_{1}$ are sums of the even and odd components ($r_{0},s_{0}\in
R_{0}$ and $r_{1},s_{1}\in R_{1}$), then $r_{0},s_{0}\in\mathrm{Z}(R)$ ensures
that the commutator%
\[
\lbrack r,s]=[r_{0}+r_{1},s_{0}+s_{1}]=[r_{1},s_{1}]=r_{1}s_{1}-s_{1}r_{1}\in
R_{0}\subseteq\mathrm{Z}(R)
\]
is a central element. Thus $R$\ is Lie nilpotent of index $2$. On the other
hand,%
\[
(r_{1})^{2}=(s_{1})^{2}=0\text{ and }0=(r_{1}+s_{1})^{2}=(r_{1})^{2}%
+r_{1}s_{1}+s_{1}r_{1}+(s_{1})^{2}%
\]
imply that $r_{1}s_{1}+s_{1}r_{1}=0$ and $[r_{1},s_{1}]=2r_{1}s_{1}$.

\bigskip

\noindent An interesting fact about the commutator parts is the following.

\bigskip

\noindent\textbf{2.5. Theorem.}\textit{ Let }$(R_{0},R_{1})$\textit{ be a
GL-grading of the }$K$\textit{-algebra }$R$\textit{ and}

\noindent$A\in\mathrm{M}_{n}(R)$\textit{ be an }$n\times n$\textit{ matrix.
Then}%
\[
AC-DA=2\lambda_{1}A_{1},
\]
\textit{where }$C,D\in\mathrm{M}_{n}([R,R])$\textit{ are the left and right
commutator parts of }$A$\textit{, }$\lambda_{1}\in R_{1}$\textit{\ is the odd
component of the symmetric determinant }$\mathrm{sdet}(A)$\textit{ and }%
$A_{1}\in\mathrm{M}_{n}(R_{1})$\textit{ is the odd component of }$A$\textit{
(}$\mathrm{sdet}(A)=\lambda=\lambda_{0}+\lambda_{1}$\textit{ and }%
$A=A_{0}+A_{1}$\textit{).}

\bigskip

\noindent\textbf{Proof.} Theorem 2.3 ensures that $AC-DA=[[\lambda,a_{i,j}]]$.
In view of our above observations concerning a GL-grading $R=R_{0}\oplus
R_{1}$, we have $[\lambda,a_{i,j}]=2\lambda_{1}a_{i,j}^{(1)}$, where
$\lambda_{1},a_{i,j}^{(1)}\in R_{1}$ are the odd components of $\lambda$ and
$a_{i,j}$ respectively. Thus $AC-DA=2\lambda_{1}A_{1}$ is proved. $\square$

\bigskip

\noindent3. THE LEFT AND RIGHT\ COMMUTATOR PARTS\ OF\ THE

\noindent\ \ \ \ CHARACTERISTIC\ MATRIX $tI_{n}-A$

\bigskip

\noindent Let $R[t]$ denote the ring of polynomials of the single commuting
indeterminate $t$, with coefficients in $R$. The symmetric characteristic
polynomial of $A$ is the symmetric determinant of the $n\times n$ matrix
$tI_{n}-A$ in $\mathrm{M}_{n}(R[t])$:%
\[
p(t)=\mathrm{tr}((tI_{n}-A)(tI_{n}-A)^{\ast})=\mathrm{tr}((tI_{n}-A)^{\ast
}(tI_{n}-A))=
\]%
\[
=\mathrm{sdet}(tI_{n}-A)=\mu_{0}+\mu_{1}t+\cdots+\mu_{n-1}t^{n-1}+n!t^{n}.
\]
The entries of the left commutator part%
\[
C(t)=n(tI_{n}-A)^{\ast}(tI_{n}-A)-\mathrm{sdet}(tI_{n}-A)I_{n}%
\]
and of the right commutator part%
\[
D(t)=n(tI_{n}-A)(tI_{n}-A)^{\ast}-\mathrm{sdet}(tI_{n}-A)I_{n}%
\]
of $tI_{n}-A$ are in the additive commutator subgroup $[R[t],R[t]]\subseteq
\lbrack R,R][t]$ of $R[t]$ generated by the elements of the form
$[f(t),g(t)]=f(t)g(t)-g(t)f(t)$ with $f(t),g(t)\in R[t]$. We also have
$\mathrm{tr}(C(t))=\mathrm{tr}(D(t))=0$.

\noindent The natural isomorphism between the rings $\mathrm{M}_{n}(R[t])$ and
$\mathrm{M}_{n}(R)[t]$ ensures that $C(t)$ and $D(t)$\ can be (uniquely)
written as (see [S3])%
\[
C(t)=C(0)+C(1)t+\cdots+C(n)t^{n}\text{ and }D(t)=D(0)+D(1)t+\cdots+D(n)t^{n},
\]
where $C(i),D(i)\in\mathrm{M}_{n}([R,R])$ and $\mathrm{tr}(C(i))=\mathrm{tr}%
(D(i))=0$\ for all $0\leq i\leq n,$

\bigskip

\noindent\textbf{Theorem} (see [S3])\textbf{.}\textit{ The following left}%
\[
(\mu_{0}I_{n}+C(0))+(\mu_{1}I_{n}+C(1))A+\cdots+(\mu_{n-1}I_{n}+C(n-1))A^{n-1}%
+(n!I_{n}+C(n))A^{n}=0
\]
\textit{and right}%
\[
(\mu_{0}I_{n}+D(0))+A(\mu_{1}I_{n}+D(1))+\cdots+A^{n-1}(\mu_{n-1}%
I_{n}+D(n-1))+A^{n}(n!I_{n}+D(n))=0
\]
\textit{Cayley-Hamilton identities hold for }$A$\textit{\ with matrix
coefficients. If }$R$\textit{\ is commutative, then}%
\[
p(t)=n!\mathrm{\det}(tI_{n}-A)\text{\textit{ , }}C(0)=C(1)=\cdots
=C(n)=D(0)=D(1)=\cdots=D(n)=0
\]
\textit{and the above identities are the }$n!$\textit{\ times scalar multiples
of the classical Cayley-Hamilton identity for }$A$\textit{.}

\bigskip

\noindent\textbf{3.1. Theorem.}\textit{ Let }$(R_{0},R_{1})$\textit{ be a
GL-grading of }$R$\textit{ and }$A\in\mathrm{M}_{n}(R)$\textit{ be an
}$n\times n$\textit{ matrix. Then }$(R_{0}[t],R_{1}[t])$\textit{ is a
GL-grading of }$R[t]$\textit{ and we have}%
\[
C(i)-D(i)-AC(i+1)+D(i+1)A=-2\mu_{i+1}^{(1)}A_{1}\text{\textit{ for all }%
}-1\leq i\leq n,
\]
\textit{where }$C(t)=C(0)+C(1)t+\cdots+C(n)t^{n}$\textit{ and }%
$D(t)=D(0)+D(1)t+\cdots+D(n)t^{n}$\textit{ are the left and right commutator
parts of }$tI_{n}-A$\textit{, }$p_{1}(t)=\mu_{0}^{(1)}+\mu_{1}^{(1)}%
t+\cdots+\mu_{n-1}^{(1)}t^{n-1}\in R_{1}[t]$\textit{\ is the odd component of
the symmetric characteristic polynomial and }$D(-1)=D(n+1)=C(-1)=C(n+1)=0,\mu
_{n}^{(1)}=\mu_{n+1}^{(1)}=0$\textit{. Notice that }$p(t)=\mathrm{sdet}%
(tI_{n}-A)=p_{0}(t)+p_{1}(t)$\textit{ and }$A=A_{0}+A_{1}$\textit{ with
}$p_{0}(t)\in R_{0}[t]$\textit{, }$A_{0}\in\mathrm{M}_{n}(R_{0})$\textit{,
}$A_{1}\in\mathrm{M}_{n}(R_{1})$\textit{.}

\bigskip

\noindent\textbf{Proof.} The fact that $(R_{0}[t],R_{1}[t])$ is a GL-grading
of $R[t]$ can be checked by straightforward calculations. We have
$(tI_{n}-A)_{0}=tI_{n}-A_{0}\in\mathrm{M}_{n}(R_{0}[t])$ and $(tI_{n}%
-A)_{1}=-A_{1}\in\mathrm{M}_{n}(R_{1}[t])$ for the even and odd components of
$tI_{n}-A$ in $\mathrm{M}_{n}(R[t])=\mathrm{M}_{n}(R_{0}[t])\oplus
\mathrm{M}_{n}(R_{1}[t])$. The application of Theorem 2.5 gives that%
\[
(tI_{n}-A)C(t)-D(t)(tI_{n}-A)=2p_{1}(t)(-A_{1})
\]
Using%
\[
(tI_{n}-A)C(t)=(tI_{n}-A)(C(0)+C(1)t+\cdots+C(n)t^{n})=
\]%
\[
=(-AC(0))+(C(0)-AC(1))t+(C(1)-AC(2))t^{2}+(C(2)-AC(3))t^{3}+\cdots
\]%
\[
\cdots+(C(n-1)-AC(n))t^{n}+C(n)t^{n+1}%
\]
and%
\[
D(t)(tI_{n}-A)=(D(0)+D(1)t+\cdots+D(n)t^{n})(tI_{n}-A)=
\]%
\[
=(-D(0)A)+(D(0)-D(1)A)t+(D(1)-D(2)A)t^{2}+(D(2)-D(3)A)t^{3}+\cdots
\]%
\[
\cdots+(D(n-1)-D(n)A)t^{n}+D(n)t^{n+1},
\]
we obtain that%
\[
(tI_{n}-A)C(t)-D(t)(tI_{n}-A)=
\]%
\[
(-AC(0)+D(0)A)+(C(0)-D(0)-AC(1)+D(1)A)t+\cdots
\]%
\[
\cdots+(C(i)-D(i)-AC(i+1)+D(i+1)A)t^{i+1}+\cdots
\]%
\[
\cdots+(C(n-1)-D(n-1)-AC(n)+D(n)A)t^{n}+(C(n)-D(n))t^{n+1}=
\]%
\[
=2p_{1}(t)(-A_{1})=-2(\mu_{0}^{(1)}A_{1}+\mu_{1}^{(1)}A_{1}t+\cdots+\mu
_{n-1}^{(1)}A_{1}t^{n-1}).
\]
Comparing the coefficients of the powers $t^{i}$, we derive that%
\[
C(i)-D(i)-AC(i+1)+D(i+1)A=-2\mu_{i+1}^{(1)}A_{1}\text{ for all }-1\leq i\leq
n.\text{ }\square
\]

\bigskip

\noindent A consequence of $(\mu_{i+1}^{(1)})^{2}=0$\ is the following.

\bigskip

\noindent\textbf{Corollary.} \textit{If }$(R_{0},R_{1})$\textit{ is a
GL-grading of }$R$\textit{ and }$A\in\mathrm{M}_{n}(R)$\textit{ is an
}$n\times n$\textit{ matrix, then the product of any two (not necessarily
different) entries of }$C(i)-D(i)-AC(i+1)+D(i+1)A$\textit{\ is zero and}%
\[
(C(i)-D(i)-AC(i+1)+D(i+1)A)^{2}=0
\]
\textit{for all }$-1\leq i\leq n$\textit{, where }$C(i),D(i)\in\mathrm{M}%
_{n}([R,R])$\textit{ are the left and right matrix coefficients in the general
Cayley-Hamilton identities.}

\bigskip

\noindent REFERENCES

\bigskip

\noindent\lbrack A] S. A. Amitsur, \textit{The T-Ideals of the Free Ring}, J.
London Math. Soc. Vol. 30 (1955), 470--475.

\noindent\lbrack D] M. Domokos, \textit{Cayley-Hamilton theorem for }%
$2\times2$\textit{\ matrices over the Grassmann algebra}, J. Pure
Appl.~Algebra 133 (1998), 69-81.

\noindent\lbrack HSvWZ] Sz. Homolya, J. Szigeti, L. vanWyk and M. Ziembowski,
\textit{Lie properties in associative algebras}, J. Algebra 573 (2021), 492--508.

\noindent\lbrack K] A. R. Kemer, \textit{Ideals of identities of associative
algebras}. Translated from the Russian by C. W. Kohls. Translations of
Math.~Monographs, 87. American Mathematical Society, Providence, RI, 1991.

\noindent\lbrack L] V. N. Latyshev,\textit{ Generalization of the Hilbert
theorem on the finiteness of bases }(Russian), Sib. Mat. Zhurn. 7 (1966),
1422-1424. Translation: Sib. Math. J. 7 (1966), 1112-1113.

\noindent\lbrack RW] D.M. Riley and M.C. Wilson. \textit{Associative rings
satisfying the Engel condition}, Proc. Amer. Math. Soc. 127 (1999), no. 4, 973-976.

\noindent\lbrack S1] J. Szigeti, \textit{New determinants and the
Cayley-Hamilton theorem for matrices over Lie nilpotent rings},
Proc.~Amer.~Math.~Soc.~125 (1997), 2245-2254.

\noindent\lbrack S2] J. Szigeti, \textit{On the characteristic polynomial of
supermatrices}, Israel J.~Math.~107 (1998), 229-235.

\noindent\lbrack S3] J. Szigeti, \textit{Cayley-Hamilton theorem for matrices
over an arbitrary ring}, Serdica Math.~J. 32 (2006), 269-276.

\noindent\lbrack S4] J. Szigeti, \textit{Integrality over fixed rings of
automorphisms in a Lie nilpotent setting}, J. Algebra 518 (2019), 198--210.

\noindent\lbrack SvW] J. Szigeti and L. van Wyk,\textit{ Determinants for
}$n\times n$\textit{ matrices and the symmetric Newton formula in the
}$3\times3$\textit{ case,} Linear Multilinear Algebra~ 62 (2014), 1076-1090.

\end{document}